\documentclass[graybox]{svmult}

\usepackage{mathptmx}       
\usepackage{helvet}         
\usepackage{courier}        
\usepackage{type1cm}        
\usepackage{makeidx}         
\usepackage{graphicx}        
\usepackage{multicol}        
\usepackage[bottom]{footmisc}
\usepackage{mathtools}

\usepackage{amsfonts}
\usepackage{amsmath}
\usepackage[utf8]{inputenc}
\usepackage[magyar]{babel}
\newcommand{\veps}{\varepsilon}

\def\be{\begin{equation}}
\def\ee{\end{equation}}

\def\bea{\begin{eqnarray}}
\def\eea{\end{eqnarray}}
\def\la{\langle}
\def\ra{\rangle}
\def\eps{\varepsilon}

\def\IP{\hbox{\rm I\kern -1.6pt{\rm P}}}
\def\IC{{\hbox{\rm
C\kern-.58em{\raise.53ex\hbox{$\scriptscriptstyle|$}}
    \kern-.55em{\raise.53ex\hbox{$\scriptscriptstyle|$}} }}}

\def\IR{\mathbb{R}}
\def\ZZ{\mathbb{Z}}
\def\IE{\mathbb{E}}
\def\IT{\hbox{\rm T\kern-.38em{\raise.415ex\hbox{$\scriptstyle|$}} }}

\def\IS{{\hbox{\rm
S\kern-.58em{\raise.53ex\hbox{$\scriptscriptstyle|$}}
    \kern-.55em{\raise.53ex\hbox{$\scriptscriptstyle|$}} }}}

\def\cR{{\cal R}}
\def\cS{{\mathcal S}}

\def\pM{\partial M}
\def\tM{\tilde M}

\def\mp{\mu^\partial}

\def\pQ{\partial Q}

  \def\length{{\rm length}}

\makeindex             

\begin{document}

\title*{Markov approximations and statistical properties of billiards}
\author{Domokos Sz\'asz}

\institute{Domokos Sz\'asz\\  Budapest University of Technology and Economics, H-1111, Egry J. u. 1,\\ \email{domaszasz@gmail.com}}
%
%

\maketitle
\today
\abstract*{}

\section{Introduction}

Mathematical billiards appeared as early as in 1912, 1913 in the works of the couple  Ehrenfest, \cite{EE12} (the wind tree model) and of D. K\"onig and A. Sz\H ucs, \cite{KSz13} (billiards in a cube) and in 1927 in the work of G. Birkhoff, \cite{B27} in (those in an oval).  Ergodic theory itself owes its birth to the desire to provide mathematical foundations to Boltzmann's celebrated  \emph{ergodic hypothesis} (cf. \cite{B31, N32, Sz11}). I briefly went over its history in my article \cite{Sz14} written on the occasion of Sinai's Abel Prize. Therefore for historic details I recommend the interested reader to consult that freely available article. Here I only mention some most relevant facts from it. In particular, I bring here two circumstances;
\begin{enumerate}
\item The two most significant problems from physics motivating the  initial study of mathematical billiards were
\begin{enumerate}
\item the ergodic hypothesis and
\item the goal to understand Brownian motion from microscopic principles.
\end{enumerate}
(In the last decades quantum billiards have also challenged both mathematicians and physicists and, moreover,  in the very last years billiard models of heat transport have also become attackable.)
\item Sinai himself was aware and highly appreciated the works of the NS Krylov, the great Russian statistical physicist who -- first of all in Russia -- brought hard ball systems, themselves hyperbolic billiards, to the attention of the community of mathematicians as a hopeful candidate for hyperbolic behavior, and possibly for ergodic one as well (cf. \cite{K79}).
\end{enumerate}
From the side of mathematics the 1960's saw the birth and rapid development of the theory of smooth hyperbolic dynamical systems with Sinai being one of the leading creators of this theory. For mathematics Sinai's 1970 paper \cite{S70} introduced a new object to study: hyperbolic billiards as hyperbolic dynamical systems with singularities.  Later it also turned out that this theory also covers basic models of {\it chaos theory}, like the Lorenz system, the H\'enon map, logistic maps, \dots.

The rich world of hyperbolic billiards and Sinai's emblematic influence on it is demonstrated by the fact that no less than three articles of this volume are devoted to Sinai's achievements in their theory. Thus I will not address here Sinai's main accomplishments in the 70's and 80's and some of their most important later expansions, which are covered in the chapter written by Leonid Bunimovich, neither will I write about the progress related to the Boltzmann-Sinai ergodic hypothesis, the topic of N\'andor Sim\'anyi's contribution. The subject of my article will be restricted to developments related to establishing {\it statistical properties of hyperbolic billiards}. These results grew out of
 \begin{itemize}
\item the appearance of the highly efficient method of Markov partitions making possible to create Markov approximations to obtain statistical properties of dynamical systems;
\item Sinai's ambition to create a mathematical theory for Brownian motion, a theory also called the dynamical theory of Brownian motion (cf. \cite{N67}). Its final goal is to derive Brownian motion from microscopic assumptions, in particular from Newtonian dynamics.
\end{itemize}

With strong simplifications our topic is the treatment of statistical properties of hyperbolic billiards via Markov approximations. The main steps in the development of this theory are, roughly speaking,  the following ones:
\begin{enumerate}
\item Markov partitions and Markov approximations for Anosov systems (and Axiom A systems) (cf. \cite{S68A,S68B,B70,S72});
\item Markov partitions and Markov approximations for 2D Sinai billiards (cf. \cite{BS80,BS81});
\item Markov sieves and Markov approximations for 2D Sinai billiards (cf. \cite{BChS90,BChS91});
\item Young's towers for hyperbolic systems with singularities, in particular for 2D Sinai billiards (cf. \cite{Y98});
\item Chernov and Dolgopyat's method of standard pairs (cf. \cite{ChD09}).
\end{enumerate}
Sinai played a founding and instrumental role in the first three steps. Chernov and Young wrote an excellent survey \cite{ChY00} on the fourth step also providing a pithy historical overview about the place of Markov partitions in the theory of dynamical systems. Referring to it permits me to focus here on the mathematical content of the exposition. At this point I note that Pesin's article \cite{P18} also in this volume discusses Markov partitions and their role in the theory of smooth hyperbolic systems in detail. My major goal in this paper will be double folded:
 \begin{itemize}
\item to put Sinai's most original achievements into perspective;
\item provide an idea about the vast and astonishing influence of them.
\end{itemize}

\section{Markov partitions for Anosov maps (and for Axiom A maps)}
Let $f$ be an Anosov diffeomorphism  of a compact differentiable manifold $M$ or an Axiom A diffeo on  $\Lambda$, one of its basic sets. Markov partitions were first constructed by Adler and Weiss \cite{AW67} (and also by by Berg \cite{B67}) for ergodic algebraic automorphisms of 2D tori. The goal of  \cite{AW67} was to provide an important positive example related to the famous isomorphism problem. Sinai's general construction for Anosov maps \cite{S68B} and its wide-ranging conclusions  \cite{S68A} revealed the sweeping perspectives of the concept. Then Bowen \cite{B70} extended the  notion to Axiom A maps and and also gave a different construction. In this section we treat both approaches simultaneously. We also remark that the content of this section finds a broader exposition in \cite{P18} in this volume.

As to fundamental notions on hyperbolic dynamical systems we refer to \cite{P18,CFS82, KH95, HK02} while here we are satisfied with a brief summary. If a diffeomorphism $f: M \to M$ has a {\it hyperbolic structure}, i. e. a decomposition into expanding vs. contracting subspaces on its unit tangent bundle $T_1 M$, then it is called an {\it Anosov map}. Then there exist two {\it foliations} into (global) stable vs. unstable invariant manifolds $\{W^u\}$ and $\{W^s\}$. Connected, bounded pieces $W^u_{\rm loc}$ (or $W^s_{\rm loc}$) of a $W^u$ (or of a $W^s$) are called {\it local stable (resp. unstable) invariant manifolds}. In particular, for any small $\veps > 0$ denote by $W^u_{(\veps)}(x)$ and $W^s_{(\veps)}(x)$  the ball-like local manifolds  of diameter $\veps$ around an $x \in W^u$. For sufficiently small $\veps$ the foliations possess a {\it local product structure}: the map $[. , .]: M \times M \to M$ is uniquely defined via $\{[x, y]\} = W^u_{(\veps)}(x)  \cap W^s_{(\veps)}(y) $. (We note that all these notions make also sense on a basic set of an Axiom A diffeo, cf. \cite{P18}.)

\begin{definition}
A subset $R$ of small diameter $\veps$ is called a {\it parallelogram} if it is closed for the operation $[., . ]$ and $R = {\rm Cl}({\rm Int} R)$.  (Further notations: $W^u_R(x) = W^u(x) \cap R$ and $W^s_R(x) = W^s(x) \cap R$.)
\end{definition}

\begin{definition}\label{Mpart}[Sinai, \cite{S68B}; Bowen, \cite{B70}]\label{def:Mpart}
A cover $\mathcal P = \{R_1, \dots, R_k\}$ of $M$ with a finite number of parallelograms with pairwise disjoint interiors is a Markov partition if for $\forall i, j$ and $\forall  \in {\rm Int}R_i \cap f^{-1}{\rm Int}R_j$ one has
\begin{enumerate}

\item $W^s_{R_i}(x) \subset f^{-1}W^s_{R_j}(fx)$
\item $W^u_{R_j}(fx) \subset fW^u_{R_i}(x)$.

\end{enumerate}
\end{definition}

The inclusions in the definition imply that, whenever $x \in {\rm Int}R_i \cap f^{-1}{\rm Int}R_j$, then $f^{-1}W^s_{R_j}(fx)$ intersects $R_i$ completely and similarly $fW^u_{R_i}(x)$ intersects ${R_j}$ completely (one can also say that these ways of intersections are Markovian).

Note that a Markov partition determines a symbolic dynamics $\tau_A$. Indeed, let $A = (a_{i, j})_{1\le i, j \le k}$ be defined  as follows: $a_{i, j} = 1\ \textrm{iff}\ {\rm Int}R_i \cap f^{-1}{\rm Int}R_j\ \neq \emptyset\ \textrm{and} \ = 0 \ \textrm{otherwise}$.  Let  $\Sigma_A$  be the subset of those $\sigma \in  \{1, \dots, k\}^\mathbb Z = \Sigma$  such that $\sigma \in \Sigma_A$ iff $\forall n \in \mathbb Z\quad a_{x_nx_{n+1}} = 1$. Then the so-called left-shift $\tau_A: \Sigma_A \to \Sigma_A$ is  defined for $\sigma \in \Sigma_A$
by $(\tau_A \sigma)_i = \sigma_{i+1}$.
$\Sigma_A$  is a closed subset of the  compact metric space $\Sigma$, a product of discrete spaces and then $\sigma_A$, called a subshift of finite type, is a homeomorphism. For a $\sigma \in \Sigma_A$ the intersection  $\cap_{i \in \mathbb Z} f^{-i}R_{\sigma_i}$ consists of a single point $x$  which we denote by $\pi(\sigma)$.
\begin{theorem} [Sinai, \cite{S68B}; Bowen, \cite{B70}]\label{thm:surjective}
$\pi: \Sigma_A \to M \ ({\textrm or}\ \Lambda)$ is a continuous surjective map and $f \circ \pi = \pi \circ \tau$.

\end{theorem}

\begin{theorem}[Sinai, \cite{S68A}]
\begin{enumerate}
\item For any transitive Anosov diffeomorphism $f$ there exists a measure $\mu^s$, positive on open subsets, such that it is invariant wrt $f$ and $f$ is a Kolmogorov-automorphism.
\item Let $\xi^s= \{W^s\}$ be the stable foliation of $M$. Then the conditional measure $\mu^s(\ .\ | W^s)$ induced on almost every $W^s$ is equivalent to the Riemannian volume on $W^s$. (Analogous statement is valid for the unstable foliation, too.)
\item If $f$ is an algebraic Anosov automorphism of $M = \mathbb T^D,\ \  D \ge 2$ (its invariant measure is Lebesgue), then $f$ is metrically conjugate (i. e. isomorphic) to a finite Markov chain.
\item The previous Markov chain has maximal entropy among all Markov chains on $\Sigma$ possessing the same possible transitions.
\end{enumerate}
\end{theorem}

Claim 1 asserts a very strong ergodic property: Kolmogorov mixing. Nevertheless, it is only a qualitative attribute, similarly to the Bernoulli property, the strongest possible ergodic one.  In the topologically mixing case an Anosov map is also Bernoulli (cf. \cite{B75}). In typical applications to problems of physics one also needs qualitative control of mixing, for instance when one has to prove a central limit theorem (CLT). In that respect Claim 3 opened principally fruitful perspectives.
Indeed, for algebraic automorphisms of $\mathbb T^D$, once they are topologically mixing, the finite Markov chain arising via the Markov partition is exponentially mixing. In such cases, if one takes a H\"older observable on $\mathbb T^D$, then this smoothness combined with the strong mixing also provides strong stochastic properties, specifically a CLT. In general, for the study of statistical properties of dynamical systems this approach makes it possible to set in the arsenal of probability theory. Later we will see the far-reaching consequences of this development. We note that in the Axiom A case, Bowen \cite{B75} established exponential correlation decay for H\"older functions and thus implying the CLT for such functions.  Claim 4 was the predecessor of Sinai's great work \cite{S72}, where by introducing symbolic dynamics in the presence of a potential function he connected the theory of dynamical systems with spin models of equilibrium statistical physics. Later this work led to thermodynamic formalism for hyperbolic systems, cf. for instance \cite{B75}.

\section{Sinai billiard and Lorentz process}

A billiard is a dynamical system describing the motion of a point particle in
a connected, compact domain $Q \subset \IT^D = \mathbb R^D / \mathbb Z^D$. In general, the boundary $\partial Q$ of
the
domain is assumed to be piecewise $C^3$-smooth; denote its smooth pieces by $\{\partial Q_\alpha| 1 \le \alpha \le J < \infty\}$. Inside $Q$ the motion is uniform while the
reflection at the boundary $\partial Q$ is elastic (by the classical rule ``the angle of incidence is
equal to the angle of reflection''). This dynamics is called the {\it billiard flow}. (In what follows we will mainly restrict our review to the discrete time billiard map.) Since the absolute value of
the velocity is a first integral of motion, the phase space of the
billiard
flow is fixed as $M=Q\times S_{D-1}$ -- in other words, every phase point
$x$ is of the form $x=(q,v)$ with $q\in Q$ and $v\in \IR^d,\ |v|=1$.
The Liouville probability measure $\mu$ on $M$ is essentially
the product of Lebesgue measures,
i.e. $d\mu= {\rm const.}\, dq dv$ (here the constant is $\frac{1}{{\textrm vol} Q\ {\textrm vol} S_{D-1}}$).

Let $n(q)$ denote the unit normal vector of a smooth component of the
boundary $\partial Q$ at the point $q$, directed inwards $Q$.
Throughout the sequel we restrict our attention to {\it dispersing
billiards}: we require that for every $q\in \partial Q$ the second fundamental
form $K(q)$ of the boundary component be positive (in fact, uniformly bounded away from $0$).

The boundary $\partial Q$ defines a natural cross-section for the billiard
flow. Consider namely
$$
\partial M = \{ (q,v) \ | \ q\in \partial Q, \ \la v,n(q)\ra \ge 0 \}.
$$
The {\it
billiard map} $T$ is defined as the first return
map on $\partial M$. The invariant measure for the map is denoted by $\mp$,
and we have $d\mp = {\rm const.} \, |\la v,n(q)\ra |\, dq dv$ (with  ${\rm const.} = \frac{2}{{\textrm vol} \partial Q\ {\textrm vol} S_{D-1}}$).
Throughout the sequel
we work with this discrete time dynamical system.

The \emph{Lorentz process} is the natural $\mathbb{Z}^D$ cover of the
above-described toric billiard. More precisely: consider $\Pi:\mathbb{R}^D
\to \mathbb{T}^D$ the factorisation by $\mathbb{Z}^D$.	Its fundamental domain
$D$ is a cube (semi-open, semi-closed) in $\mathbb{R}^D$, so $\mathbb{R}^D =
\cup_{z \in \mathbb{Z}^D} (D+z)$, where $D+z$ is the translated fundamental
domain.
We also lift the scatterers to $\mathbb R^D$ and define the phase space of the Lorentz flow as $\tilde M = \tilde Q\times S_{D-1}$, where $\tilde Q = \cup_{z \in \mathbb{Z}^D} (Q+z)$. In the non-compact space $\tilde M$ the dynamics is denoted by  $\tilde S^t$  and the billiard map on $\partial \tilde M$  by	$\tilde T$; their natural projections to the configuration space $\tilde Q$ are denoted by $L(t) = L(t;x), t \in \IR_+$ and $L^\partial_{n} \in \ZZ_+$  and called (periodic) Lorentz flows or processes with natural invariant measures $\tilde \mu$ and $\tilde \mu^\partial$, respectively.

The \emph{free flight
  vector} $\tilde\psi: \tilde M \to \mathbb{R}^D$ is defined as follows:
$\tilde \psi(\tilde x)=\tilde q(\tilde T\tilde x)-\tilde q (\tilde x)$.

\begin{definition}
  The Sinai-billiard (or the Lorentz process) is said to have \emph{finite horizon} if the free flight vector
  is bounded.  Otherwise the system is said to have \emph{infinite horizon}.
\end{definition}

\subsection{Singularities.}

\subsubsection{Tangential singularities} Consider the set of tangential reflections, i.e.
$$
\cS := \{ (q,v)\in \partial M \ |\  \la v, n(q)\ra =0 \}.
$$
It is easy to see that the map $T$ is not continuous at the set
$T^{-1}\cS$. As a consequence, the (tangential) part of the singularity set for iterates
$T^n,\ n \ge 1$ is
$$
\cS^{(n)}= \cup_{i=1}^n \cS^{-i},
$$
where in general $\cS^k=T^k \cS$.

\subsubsection{Multiple collisions}
After the billiard trajectory hits $\partial Q_{\alpha_1} \cap \partial Q_{\alpha_2}$ (for some $\alpha_1 \neq \alpha_2$), the orbit stops to be uniquely defined and there arise two - or more - trajectory branches. Denote
$$
\cR := \{ (q,v)\in \partial M \ |\  q \in \partial Q_{\alpha_1} \cap \partial Q_{\alpha_2}\ \textrm{for some}\  \alpha_1 \neq \alpha_2  \}.
$$
{\bf Standing assumption.} We always assume that if  $q \in \partial Q_{\alpha_1} \cap \partial Q_{\alpha_2}\ \textrm{for some}\  \alpha_1 \neq \alpha_2$, then these two smooth pieces meet in $q$ in general position (in the planar case this implies a non-zero angle between the pieces).

It is easy to see that the map $T$ is not continuous at the set
$T^{-1}\cR$. As a consequence, for iterates
$T^n,\ n \ge 1$ the part of the singularity set, caused by multiple collisions, is
$$
\cR^{(n)}= \cup_{i=1}^n \cR^{-i},
$$
where in general $\cR^k=T^k \cR$.

\subsubsection{Handling the singularities}
Here we only give a very rough idea. When hitting any type of singularities the map is not continuous (the flow is still continuous at tangential collisions, but it stops being smooth). Consequently, $W^u_{\rm loc}$ are those connected pieces of $W^u$ which never hit
\[
\Sigma_{n=1}^\infty (\cS^{(n)} \cup \cR^{(n)})
\]
in the future. (Reversing time one obtains $W^s_{\rm loc}$). The basic observation in Sinai's approach to billiards was that these smooth pieces are $D-1$-dimensional local manifolds almost everywhere.

An additional difficulty connected to tangential singularities is that the expansion rate in the direction orthogonal to the singularities is infinite, a phenomenon breaking the necessary technical quantitative bounds. The way out was found in \cite{BChS91} where the authors introduced additional -- so called secondary -- singularities. These will further cut $W^{u, s}_{\rm loc}$ and in what follows  $W^{u, s}_{\rm loc}$ will denote these smaller pieces, themselves $D-1$-dimensional local manifolds almost everywhere. (Detailed exposition of these can be found in \cite{ChM06}  in the planar case, and  in \cite{BChSzT03} in the multidimensional case.)

\section{Statistical properties of 2D periodic Lorentz processes}\label{sec:StatProp}
Given the successes of Markov partitions for smooth hyperbolic systems and of Sinai's theory of ergodicity for hyperbolic billiards, a prototype of hyperbolic systems with singularities, it is a natural idea to extend the method of Markov partitions to Sinai billiards. Yet, when doing so there arise substantial difficulties. The most serious one is that basic tools of hyperbolic theory: properties of the holonomy map (also called canonical isomorphism), distorsion bounds, etc. are only valid for smooth pieces of the invariant manifolds (maximal such components are called local invariant manifolds and denoted by $W^{u, s}_{\textrm loc}$, cf. 3.1.3). These can, however, be arbitrarily short implying that a Markov partition can only be infinite, not finite. Technically this means the construction of a countable set of parallelograms (products of Cantor sets in this case) with an appropriate Markov interlocking; this is a property which is really hard to control.

Assume we are given a Sinai billiard. In the definition of parallelogram we make two changes. First, in the operation $[. , . ]$ we only permit $W^{u, s}_{\textrm loc}$ and, second, we do not require $R = {\rm Cl}({\rm Int} R)$ any more. Now we will denote $W^{u,s}_R(x) = W^{u, s}_{\textrm loc}(x) \cap R$.
\begin{definition}\label{Mpart}[Bunimovich-Sinai, \cite{BS80}, Bunimovich-Chernov-Sinai, \cite{BChS90}]

A cover $\mathcal P = \{R_1, R_2, \dots, \}$ of $M$ with a countable number of parallelograms, satisfying $\mu_1(R_i \cap R_j) = 0$  ( $ \forall 1 \le  i < j$), is a Markov partition if  one has that  if $x \in {\rm Int}R_i \cap f^{-1}{\rm Int}R_j$, then
\begin{enumerate}

\item $W^s_{R_i}(x) \subset f^{-1}W^s_{R_j}(fx)$
\item $W^u_{R_j}(fx) \subset fW^u_{R_i}(x)$.

\end{enumerate}

\end{definition}

From now on we assume that $D = 2$ and that, unless otherwise stated, the {\it horizon is finite}.

\subsection{Bunimovich-Sinai, 1980}
\begin{theorem}[Bunimovich-Sinai, \cite{BS80}]\label{thm:MP80}
Assume  that for the billiard in $Q = \mathbb T^2 \setminus \Sigma_{j=1}^J \mathcal O_j$ the strictly convex obstacles $\mathcal O_j$ are closed, disjoint with  $C^3$-smooth boundaries. Then for the billiard map $T$ there exists a countable Markov partition of arbitrarily small diameter.
\end{theorem}
It is worth mentioning that the statement of Theorem \ref{thm:surjective} still keeps holding for the constructed Markov partition. (We also note that a correction and a simplification of the construction was given in \cite{BS86,L86}.)\\
In the companion paper \cite{BS81} the authors elaborate on further important properties of the constructed Markov partition and prove groundbreaking consequences for the Lorentz process. For $x=(q, v) \in \pM$ denote $T^nx = (q(n), v(n))$ and for $x \in \tM$ the - diffusively - rescaled version of the Lorentz process by
\[
L_A(x) = \frac{1}{\sqrt A} L(At; x)\qquad (t \in \IR_+).
\]

\begin{theorem}[Bunimovich-Sinai, \cite{BS81}]\label{thm:CD80}
There exists a constant $\gamma \in (0, 1)$ such that for all sufficiently large $n$
\[
\left|\mathbb  E_{\mp}(v(0) v(n))\right| \le \exp(-n^\gamma).
\]

\end{theorem}
The proof uses Markov approximation. One of its essential elements is that a rank function is introduced on elements of the partition: roughly speaking the smaller the element is the larger is its rank. Though the tail distribution for the rank is exponentially decaying nevertheless the well-known Doeblin condition of probability theory, ensuring exponential relaxation to equilibrium, does not hold for one step transition  probabilities. Fortunately it does hold for higher step ones, still with the step size depending on the rank of the element of the partition.
This weaker form of Doeblin property implies that $\gamma$ is necessarily smaller than $1$.
Yet this is a sufficiently strong decay of correlations to imply convergence to Brownian motion. Assume $\nu$ is a probability measure on $\tilde  M$ supported on a bounded domain and  absolutely continuous wrt $\tilde \mu$.
\begin{theorem}[Bunimovich-Sinai, \cite{BS81}]\label{thm:BM80}
With respect to the initial measure $\nu$, as $A \to \infty$
\[
L_A(t) \Rightarrow B_\Sigma (t)
\]
where $B_\Sigma (t)$ is the planar Wiener process with zero shift and covariance matrix $\Sigma$ and the convergence is weak convergence of measures in $C[0, 1]$ (or in $C[0, \infty]$).\\
Moreover, if the Lorentz process is not localised and the scatterer configuration is symmetric wrt the line $q_x = q_y$, then $\Sigma$ is not singular.
\end{theorem}

\subsection{Bunimovich-Chernov-Sinai, 1990-91}
Ten years after the first construction Bunimovich, Chernov and Sinai revisited the topic in two companion papers. The authors not only simplified the original constructions and proofs of \cite{BS80,BS81}, but also clarified and significantly weakened the conditions imposed. Below we summarize the most important attainments. \\
{\bf Wider class of billiards.} Consider a planar billiard in $Q \subset \IT$ with piecewise $C^3$-smooth boundary. Impose the following conditions:
\begin{enumerate}
\item If $q \in \partial Q_{\alpha_1} \cap \partial Q_{\alpha_2}$ for some $\alpha_1 \neq \alpha_2$, then the angle between $\partial Q_{\alpha_1}$ and  $\partial Q_{\alpha_2}$ is not zero;
\item There exists a constant $K_0 = K_0(Q)$ such that the multiplicity of the number of curves of $\cup_{i=-n}^n (\cS^{-i} \cup \cR^{-i})$ meeting at any point of $\pM$ is at most $K_0n$;
\end{enumerate}

\begin{theorem}[Bunimovich-Chernov-Sinai, \cite{BChS90}]\label{thm:MP90}
Assume  that a planar billiard with finite horizon satisfies the two conditions above. Then for the billiard map $T$ there exists a countable Markov partition of arbitrarily small diameter.
\end{theorem}

Extending Theorems \ref{thm:CD80} and \ref{thm:BM80} and simplifying the proofs the three authors could get the following results.
Denote by $\mathcal H_\beta, \ \ \beta > 0$ the class of $\beta$-H\"older functions $h: \pM \to \IR$ (i. e. $\exists C(h)$ such that $\forall \alpha \le J\ \ {\textrm and}\ \ \forall x,y \in \pQ_\alpha\ \ \  |h(x) - h(y)| \le C(h)|x-y|^\beta$). We note that the sequence $X_n = h(T^nx)\ (n \in \ZZ)$ is stationary wrt the invariant meassure $\mp$ on $\pM$.

\begin{theorem}[Bunimovich-Chernov-Sinai, \cite{BChS91}]\label{thm:CD91}
Assume the billiard satisfies the previous two conditions and take a function $h \in \mathcal H_\beta$ with $\mathbb E_{\mp} h = 0$. Then $\forall n \in \ZZ$ one has
\[
|\mathbb E_{\mp} X_0X_n| \le C(h) e^{-a\sqrt n}
\]
where $a = a(Q) > 0$ only depends on the billiard table.
\end{theorem}

\begin{theorem}[Bunimovich-Chernov-Sinai, \cite{BChS91}]\label{thm:BM91}
Assume the billiard satisfies the previous two conditions.
With respect to the initial measure $\nu$, as $A \to \infty$
\[
L_A(t) "\Rightarrow" B_\Sigma (t)
\]
where $B_\Sigma (t)$ is the planar Wiener process with zero shift and covariance matrix $\Sigma$ and the convergence $"\Rightarrow"$ is weak convergence of finite dimensional distributions.
Moreover, if the Lorentz process is not localised, then $\Sigma$ is not singular.
\end{theorem}
It is important to add that these two papers also discuss semi-dispersing billiards and, in general, provide a lot of important information about the delicate geometry of various examples. The proofs of Theorems \ref{thm:CD91} and \ref{thm:BM91} do not use the Markov partition of \ref{thm:MP90} directly but build up a Markov approximation scheme by using so-called approximate finite Markov-sieves. An immediate additional success of the method of Markov sieves was a spectacular physical application. In \cite{ChELS93,ChELS93B} the authors could study a billiard-like models under the simultaneous action a Gaussian thermostat and a small external field. The interesting feature of the model is that the system is not Hamiltonian and has an attractor. Among other beautiful results they derive a formula for the rate of entropy production and verify Einstein's formula for the diffusion coefficient.

\section{Further progress of Markov methods}

\subsection{Markov towers} Sinai's ideas on connecting dynamical systems with statistical physics and probability theory have been justified by the works mentioned in Section \ref{sec:StatProp}. It became evident that billiard models are much appropriate for understanding classical questions of statistical physics. Having worked out and having simplified the meticulous details of Markov approximations the way got opened for further progress.

 Young, who had also had experience with other hyperbolic systems with singularities, like logistic maps and the H\'enon map, was able to extract the common roots of the models. She introduced a fruitful and successful system of axioms under which one can construct Markov towers, themselves possessing a Markov partition. A major advantage of her approach was the following:  the papers discussed in Section \ref{sec:StatProp} had showed that, though it was indeed possible, but at the same time rather hard to construct Markov partitions for billiards directly. An important idea of \cite{Y98} is that one can rather use renewal properties of the systems and build Markov towers instead. A remarkable accomplishment of the tower method was that Young could improve the stretched exponential bound on correlation decay of Theorem \ref{thm:CD91} to the optimal, exponential one. Since this was just one - though much important - from the applications of her method, in \cite{Y98} she restricted her discussion to the case of planar finite-horizon Sinai billiards with $C^3$-smooth scatterers.

 \begin{theorem}[Young, \cite{Y98}]\label{thm:Y98}
Assume the conditions of Theorem \ref{thm:MP80}. Then for any $\beta > 0$ and for any $g, h \in \mathcal H_\beta$ there exist $a > 0$ and $C = C(g, h)$ such that
\begin{enumerate}
\item
$\forall n \in \ZZ$ one has
\[
|\mathbb E_{\mp} (g \circ f^n) h - E_{\mp} (g) E_{\mp} (h)| \le C e^{-an};
\]
\item
\[
\frac{1}{\sqrt n}\left(\sum_0^{n-1}g \circ f^n - n E_{\mp} (g)\right)  \xrightarrow{distr} \mathcal N(0, \sigma)\qquad (n \to \infty)
\]
where $N(0, \sigma)$ denotes the normal distribution with $0$ mean and variance $\sigma^2 \ge 0$.
\end{enumerate}
\end{theorem}

 \begin{remark}
Pay attention to the differences in the assertions of  Theorems \ref{thm:CD80}, \ref{thm:BM91} and \ref{thm:Y98}. The authors of the last two works have not claimed and checked tightness which was, indeed, settled in \cite{BS81}.
\end{remark}

 Beside the original paper one can also read \cite{ChY00} describing the ideas in a very clear way. It is also worth noting that \cite{Y99} extended the tower method to systems where the renewal time has a tail decreasing slower than exponential. Young's tower method can be considered as a fulfillment of Sinai's program. Her axioms for hyperbolic systems with singularities serve as an autonomous - and most popular - subject and make it possible to discuss wide-ranging delicate stochastic properties of the systems covered, interesting either from probabilistic or dynamical or physical point of view. Two examples from the numerous applications are \cite{R-BY08} proving large deviation theorems for systems satisfying Young's axioms and \cite{PS10} describing a recurrence type result in the planar Lorentz process setup.

Unfortunately without further assumptions the method works so far  for the planar case, only.  However, in their paper \cite{BT08} providing an important achievement, B\'alint and T\'oth formulated a version of the tower method for multidimensional billiards under the additional 'complexity' hypothesis, whose verification for multidimensional models is a central outstanding question of the theory.

\subsection{Standard pairs} Another astonishing development of Markovian tools was the 'standard pair' method of Chernov and Dolgopyat \cite{ChD09}. This method has already had remarkable applications but so far it is not easy to see where its limits are. As to a recent utilization we can, for instance, mention that standard pairs have also been applied to the construction of SRB measures for smooth  hyperbolic maps in any dimension; section 3 of \cite{CDP14} provides a brief introduction to the tool, too. Since -- for systems with singularities -- until now  the method of standard pairs does not have a clear survey exposition as  \cite{ChY00} is for the tower method, very briefly we present a theorem showing how it handles Markovity.

Let $(\pM, T, \mp)$ be the billiard ball map - for simplicity for a planar billiard. A standard pair is $\ell=(W, \rho)$ where $W$ is an unstable curve, $\rho$ is a nice probability density on $W$ (an unstable curve is a smooth curve in $\pM$ whose derivatives at every point lie in the unstable cone).
Decompose $\pM$  into a  family of nice standard pairs. Select a  standard pair  $\ell=(W, \rho)$ from this family. Fix a nice function $A: M \to \mathbb R$. Then according to the well-known law of total probability
\begin{equation}\label{eq:partitioning}
 \IE_\ell(A\circ T^n)=\sum_\alpha c_{\alpha n}
\IE_{\ell_{\alpha n}} (A)
\end{equation}
where $c_{\alpha n}> 0$,
$\sum_\alpha c_{\alpha n}=1$. The $T^n$-image of $W$ is cut to a finite or countable number of pieces $W_{\alpha n}$. Thus
$\ell_{\alpha n}=(W_{\alpha n}, \rho_{\alpha n})$
are disjoint standard pairs with $T^n W = \cup_\alpha W_{\alpha, n}$ where
 $\rho_{\alpha n}$ is the pushforward of $\rho$ up
to a multiplicative factor.

\begin{theorem}[Chernov-Dolgopyat \cite{ChD09}; Growth lemma\ $\sim$\ Markov property]
If $n \ge \beta_3 |\log\length(\ell)|$, then
$$\sum_{\length(\ell_{\alpha n})<\eps} c_{\alpha n}
\leq \beta_4 \eps . $$
\end{theorem}
Equation \ref{eq:partitioning} expressed how an unstable curve is partitioned after $n$ iterations. Among the arising pieces there are, of course,longer and shorter pieces. The theorem provides a quantitative estimate for the total weight of pieces shorter than $\veps$.
This theorem, a quantitative formulation of Sinai's traditional billiard philosophy: 'expansion prevails partitioning' was not new, in various forms it had appeared in earlier works, too. Its consequent application, however, together with modern formulations of averaging theory and a perturbative study of dynamical systems was absolutely innovative and most successful.

\section{Further successes  of Markov methods}
Because of the abundance of related results my summary will be very much selective. My main guiding principle will be that I try to focus on those developments that are either directly related to Sinai's interests or even to the problems he raised or alternatively show a variety of questions from physics.
\subsection{Applications of the tower method}
\begin{enumerate}
\item In 1999 already, Chernov \cite{Ch99} could extend the exponential correlation bound of \cite{Y98} to planar billiards with infinite horizon: it holds for H\"older observables. (The work also contains precious analysis of the growth lemma, of homogeneity layers, ...).
\item As mentioned before the works \cite{ChELS93,ChELS93B} treated billiards with small external forces. Chernov, in a series of articles (that started with \cite{Ch01} and ended with \cite{ChZhZh13}) worked out a comprehensive  theory of these models.
\item The methods initiated by Sinai also made it possible to study Sinai billiards with small holes, a model suggested by physicists. Early answers to the questions were treated by Markov partitions (cf. \cite{ChM97,ChMT00} in case of Anosov maps) whereas later, results for billiards were found by applying Young towers (e. g. in \cite{DWY10,DWY12}).

\item After the CLT of \cite{BS81} for the Lorentz process, Sinai formulated the question: is P\'olya recurrence true for finite horizon planar billiards? Positive answers were obtained by \cite{S98}, \cite{C99} and \cite{SzV04}. The latter work accomplished that by proving a local version of the CLT for the planar finite horizon Sinai billiard (cf. next point, too).
\item In the planar infinite horizon case the free flight function determining the Lorentz process is not H\"older (not bounded even), so the correlation bound of \cite{Ch99} is not applicable to it. In fact, in this case the scaling of the Lorentz process, in a limit law like that of Theorems \ref{thm:BM80} and \ref{thm:BM91}, is different and it is  $\sqrt{n \log n}$ rather than $\sqrt{n }$. This was shown in  \cite{SzV07}, where by extending  the method of \cite{SzV04} P\'olya recurrence was also obtained via an appropriate local limit theorem. The analogous  results for the Lorentz flow with many other interesting theorems -- also in the presence of external field -- were obtained in \cite{ChD09B}.
    \item As this was observed in \cite{BG06}, in limit laws for stadium billiards there may arise limit theorems with both classical and non-classical scaling (cf. previous point).
\end{enumerate}

\subsubsection{Applications of the method of standard pairs}
Roughly saying the additional step, this method represents in the evolution of  Markov approximation tools, can be characterized by two main intertwining advantages: first it makes possible to treat systems with  two (or several) time scales and second it is appropriate for a perturbative description of dynamical systems, in particular, of billiards. The basic reference is \cite{ChD09} though the authors started lecturing about it as early as in 2005 (see, for instance, \cite{D05}).

\begin{enumerate}
\item \cite{ChD09} makes an important step in the dynamical theory of Brownian motion: two particles move on a planar Sinai billiard table, with one of them: an elastic disk being much heavier than the other one: a point particle. Since the motion of the heavy disk is slow, for the point particle -- in short time intervals -- statistical properties hold (among the scatterers of the original Sinai billiard plus the -- temporarily fixed -- heavy disk particle). However, when the heavy particle gets close to any of  the original scatterers, then additional phenomena appear and so far this is the limit of the applicability of the method. \cite{BChD11} is the first step toward extending the time interval where the theory is hoped to be applicable.
\item Multiple times scales are treated by standard pairs in \cite{DL11}. Though their model is not a billiard one, actually the dynamics is smooth, nevertheless the work is much successful in deriving a mesoscopic, stochastic process from Newtonian, microscopic laws of motion. This task, also important in a rigorous study of a heat transport model of physicists (cf.   \cite{GG08}), is the subject of \cite{BNSzT15} for a billiard model.
\item Another spectacular development was obtained in \cite{ChD09C}. The authors were considering a point particle moving in ${\IR}^2$ in the presence of a constant force among periodically situated strictly convex scatterers (the horizon is assumed to be finite). They could derive non-classical limit laws  both  for  the velocity and the position of the particle and moreover, they could also prove the recurrence of the particle.
    \item Sinai raised the following problem in 1981: consider a finite horizon, planar Sinai billiard and displace one scatterer a bit. Prove for it an analogue of Theorem \ref{thm:BM91}. This problem was answered by the method of standard pairs in the companion works \cite{DSzV08,DSzV09}.
    \item Returning to heat conduction: in \cite{DN15} the authors could derive the heat equation for a Lorentz process in a quasi one-dimensional tube being long finite and asymptotically infinite. The boundaries, on the one hand, absorb  particles reaching them and, on the other hand, particles are also injected with energies corresponding to different temperatures.
\end{enumerate}

Closing this section I note that despite the striking successes of Markov approximations techniques their applicability so far is essentially restricted to two dimensional models, except perhaps for \cite{CDP14}. Sinai's original works addressed explicitly planar models, only. Though the tower method is extended to the multidimensional case, the 'complexity condition' arising in it has not been checked hitherto for any multidimensional system yet it is strongly believed that it does hold at least typically.

As mentioned earlier many successes of the theory of hyperbolic billiards were motivated by problems of physics. The recent survey \cite{D14}, however, shows  that, as usual, there are more open problems than those solved.

{\bf Acknowledgement.} The author is highly indebted to Yasha Pesin for his careful reading of the manuscript and many useful remarks.

\end{document}